\documentclass[11pt, reqno]{amsart}
\usepackage{amsmath}
\usepackage{amssymb}

\font\tengothic=eufm10
\font\sevengothic=eufm7
\newfam\gothicfam
      \textfont\gothicfam=\tengothic
      \scriptfont\gothicfam=\sevengothic
\def\goth#1{{\fam\gothicfam #1}}

\numberwithin{equation}{section}

\newenvironment{mydef}{{\bf Definition.}}{\hspace*{\fill} \par\vspace{1ex}}

\begin{document}
\setlength{\baselineskip}{1.7em}
\newtheorem{thm}{\bf Theorem}[section]
\newtheorem{pro}[thm]{\bf Proposition}
\newtheorem{claim}[thm]{\bf Claim}
\newtheorem{lemma}[thm]{\bf Lemma}
\newtheorem{cor}{\bf Corollary}[thm]
\newtheorem{remark}[thm]{\bf Remark}
\newtheorem{example}[thm]{\bf Example}
\newcommand{\cc}{{\mathbb C}}
\newcommand{\pp}{{\mathbb P}}
\newcommand{\zz}{{\mathbb Z}}
\newcommand{\nn}{{\mathbb N}}
\newcommand{\R}{{\mathcal R}}
\newcommand{\A}{{\mathcal A}}
\newcommand{\x}{{\mathbb X}}
\newcommand{\y}{{\mathbb Y}}
\newcommand{\ix}{I_{\mathbb X}}
\newcommand{\ci}{{\mathcal I}}
\newcommand{\B}{{\bf B}}
\newcommand{\bi}{{\bf I}}
\newcommand{\V}{{\bf V}}
\newcommand{\cv}{{\mathcal V}}
\newcommand{\bL}{{\bf L}}
\newcommand{\co}{{\mathcal O}}
\newcommand{\cf}{{\mathcal F}}
\newcommand{\cl}{{\mathcal L}}
\newcommand{\cm}{{\mathcal M}}
\newcommand{\cn}{{\mathcal N}}
\newcommand{\spec}{\mbox{Spec }}
\newcommand{\pr}{\mbox{Proj }}
\newcommand{\sv}{S_{V_\lambda}}
\newcommand{\m}{--\!\!\rightarrow}
\newcommand{\smap}{\rightarrow\!\!\!\!\!\rightarrow}
\newcommand{\sfrac}[2]{\frac{\displaystyle #1}{\displaystyle #2}}
\newcommand{\under}[1]{\underline{#1}}
\newcommand{\ov}[1]{\overline{#1}}
\newcommand{\sh}[1]{{\mathcal #1}}
\newcommand{\lex}{{\le}_{\mbox{\scriptsize lex}}}

\title{Arithmetic Macaulayfication of Projective Schemes}
\author{Steven Dale Cutkosky and Huy T\`ai H\`a}
\thanks{Steven Dale Cutkosky was partially supported by NSF}
\thanks{{\it 2000 Mathematical Subject Classification}. 14E05, 13A30}
\keywords{Arithmetic Macaulayfication, Rees algebra, blowing up}
\address{Department of Mathematics, University of Missouri, Columbia MO 65201, USA}
\email{dale@math.missouri.edu}
\email{tai@math.missouri.edu}
\begin{abstract}
In this paper, we study  arithmetic Macaulayfication of projective schemes and Rees
algebras of ideals. We discuss the existence of an arithmetic Macaulayfication for
projective schemes. We give a simple neccesary and sufficient condition for
nonsingular  projective varieties to possess an arithmetic Macaulayfication
(Theorem \ref{ACM}). We also show that this condition is  sufficient in general,
but give examples to show that it is not in general necessary. We further consider
Rees algebras $\R_\lambda(I) = R[(I_{\lambda})t]$ (truncated Rees algebras)
associated to a homogeneous ideal $I$ and show that they are Cohen-Macaulay for
large $\lambda$ in some important cases (Theorem \ref{trunc-rees} and Corollary
\ref{CorollaryLCI}).
\end{abstract}
\maketitle

\begin{center}
{\it Dedicated to Wolmer Vasconcelos on the occation of his sixty fifth birthday.}
\end{center}

\setcounter{section}{-1}
\section{Introduction}
In this paper, we study  arithmetic Macaulayfication of
projective schemes and Rees algebras of ideals.

In the first part of the paper, we discuss the problem of arithmetic
Macaulayfication of projective schemes. This is a globalization of the problem of
arithmetic Macaulayfication of local rings, which was first considered by Barshay
in \cite{barshay}, and then studied extensively by many authors, such as Goto and
Shimoda \cite{gs}, Goto and Yamagishi \cite{gy}, Brodmann \cite{br}, Schenzel
\cite{sch}, Lipman \cite{l}, Aberbach \cite{ab}, Kurano \cite{ku}, Aberbach, Huneke
and Smith \cite{ahs}, and finally solved by Kawasaki \cite{ka}. We give a neccesary
and sufficient condition for a nonsingular projective scheme over a field $k$ of
characteristic 0 to have an arithmetic Macaulayfication.

\begin{thm} (Theorem \ref{nonsingular})
Suppose $X$ is a nonsingular projective scheme over a field $k$ of characteristic 0.
Then, $X$ has an arithmetic Macaulayfication if and only if $H^0(X,\co_X)=k$ and
$H^i(X, \co_X) = 0$ for all $i =1, \ldots, \dim X-1$.
\end{thm}

We show that the cohomological conditions of Theorem \ref{nonsingular} are
sufficient conditions for an unmixed projective scheme to have an arithmetic
Macaulayfication (Theorem \ref{ACM}). This result follows from the work of Kawasaki
(\cite{ka1}, \cite{ka}). However, we show that the cohomological conditions of
Theorem \ref{nonsingular} are   not necessary in general (Example \ref{Example1},
Example \ref{ExampleB}).

In the second part of this paper, we consider a natural class of Rees algebras
associated to a homogeneous ideal, the {\it truncated Rees algebras}. This class of
Rees algebras was first considered by the second author in \cite{ha-thesis} and
\cite{ha-rees} for the defining ideal of a set of points in $\pp^2$. It gave a new
tool to completely answer the question on defining equations of projective
embeddings of certain rational surfaces (see \cite[Section 4.3]{ha-thesis}).

Our main result of this section is Theorem \ref{trunc-rees}, from which we can
conclude results such as the following.

\begin{cor}(Corollary \ref{CorollaryLCI}) Suppose that $X$ is a projective
Cohen-Macaulay scheme over a field $k$ such that $H^0(X,\co_X)=k$,
$H^i(X,\co_X)=0$ for $i>0$ and $\sh{I}\subset \co_X$ is an ideal sheaf which is locally
a complete intersection. Then
\begin{enumerate}
\item There exists a Cohen-Macaulay standard graded $k$-algebra $R$ with $a(R)<0$
such that $X\cong \pr R$.
\item If $I\subset R$ is a homogeneous ideal such that $\tilde I\cong \sh{I}$,
then there exists $\lambda_0\ge\delta(I)+1$ (where $\delta(I)$ is the maximum degree
of a minimal set of homogeneous generators of $I$) such that the truncated Rees
algebra $$ \R_\lambda(I) = R[(I_\lambda) t] $$ is Cohen-Macaulay for all
$\lambda\ge\lambda_0$.
\end{enumerate}
\end{cor}

To prove Theorem \ref{trunc-rees}, we combine the method of \cite{hyry} for
studying the local cohomology of multigraded algebras with the results of \cite{ch}.

Throughout this paper, let $k$ be a field. We follow the notations of \cite{ei} and
\cite{hart}.

\section{Arithmetic Macaulayfication}

Suppose that $X$ is a projective scheme over a field $k$. We will say that $X$ is
{\it arithmetically Cohen-Macaulay} if there exists a Cohen-Macaulay standard graded
$k$ algebra $S$ such that $X\cong \pr S$.

\begin{mydef} Suppose that $X$ is a projective scheme over a field $k$.
An {\it arithmetic
Macaulayfication} of $X$ is a proper
birational morphism $\pi: Y \rightarrow X$ such that $Y$ is arithmetically
Cohen-Macaulay.
\end{mydef}

We shall first prove a very basic result on arithmetically Cohen-Macaulay schemes.

\begin{lemma} \label{cond}
\begin{enumerate}
\item Suppose that $Y = \pr S$ is an arithmetically Cohen-Macaulay scheme. Then, $Y$ is a
Cohen-Macaulay scheme,  $H^i(Y, \co_Y) = 0$ for $i = 1, \ldots, \dim Y-1$,
and $H^0(Y,\co_Y)=k$.
\item Suppose that $Y = \pr S$ is a Cohen-Macaulay scheme, $H^i(Y, \co_Y) = 0$ for
$i =1, \ldots, \dim Y-1$ and $H^0(Y,\co_Y)=k$. Then, there exists an integer $n_0$ such that for all $n
\ge n_0$, the Veronese embedding of $Y$ by $H^0(Y,\co_Y(n))$ is arithmetically
Cohen-Macaulay.
\end{enumerate}
\end{lemma}

\begin{proof} Let $\goth{m}$ be the maximal homogeneous ideal of $S$.
We have isomorphisms
\[ \oplus_{n \in \zz} H^i(Y, \co_Y(n)) \cong H^{i+1}_\goth{m}(S), \forall \ i \ge 1, \]
and an exact sequence
\[ 0 \rightarrow H^0_\goth{m}(S) \rightarrow S \rightarrow \oplus_{n \in \zz} H^0(Y,
\co_Y(n)) \rightarrow H^1_\goth{m}(S) \rightarrow 0. \]
(1) is immediate since $S$ is Cohen-Macaulay if and only if $H^i_\goth{m}(S)=0$ for
$i\le d=\text{dim }\pr S$.

To prove (2) we first observe that $H^i(Y,\co_Y(n))=0$ for $i>0$ and $n \gg 0$ by
Serre vanishing. Since $Y$ is Cohen-Macaulay, we also have $H^i(Y,\co_Y(n))=0$ for
$i<d$ and $n \ll 0$.
\end{proof}

\begin{thm} \label{hiro}(Hironaka \cite{h})
Suppose that $\pi: Y \rightarrow X$ is a birational morphism of projective
nonsingular varieties over a field $k$ of characteristic 0. Then
\[ H^i(Y, \co_Y) \simeq H^i(X, \co_X) \ \forall \ i. \]
\end{thm}

\begin{proof} By  resolution of indeterminancy (\cite{h}), there
exists a commutative diagram of projective morphisms
\[ \begin{array}{rrl} Z & \stackrel{f}{\rightarrow} & Y \\
& g \searrow & \downarrow \pi \\
& & X \end{array}
\]
such that $g$ is a product of blowups of nonsingular subvarieties,
\[ g: Z = Z_n \stackrel{g_n}{\rightarrow} Z_{n-1} \stackrel{g_{n-1}}{\rightarrow}
\ldots \stackrel{g_2}{\rightarrow} Z_1 \stackrel{g_1}{\rightarrow} Z_0 = X. \]
We have (\cite{mat} or Lemma 2.1 \cite{ch})
\[ R^i {g_j}_{*} \co_{Z_j} = \left\{ \begin{array}{ll}
0, & i > 0 \\ \co_{Z_{j-1}}, & i = 0 \end{array} \right. \]
Thus,
\begin{eqnarray}
R^ig_{*} \co_Z & = & \left\{ \begin{array}{ll} 0, & \mbox{ if } i > 0 \\ \co_X, &
\mbox{ if } i = 0 \end{array} \right. \label{d2}
\end{eqnarray}
and $$ g^*:H^i(X, \co_X) \cong H^i(Z, \co_Z) $$ for all $i$. Now, by considering
the commutative diagram
\[ \begin{array}{rrc} H^i(Z, \co_Z) & \stackrel{f^{*}}{\leftarrow} & H^i(Y, \co_Y) \\
& g^{*} \nwarrow & \uparrow \pi^{*} \\
& & H^i(X, \co_X) \end{array}
\]
we conclude that $\pi^{*}$ is one-to-one. To show that $\pi^{*}$ is an isomorphism
we now only need to show that $f^{*}$ is also one-to-one.

By resolution of indeterminancy it also gives a new diagram
\[ \begin{array}{rrl} W & \stackrel{\gamma}{\rightarrow} & Z \\
& \beta \searrow & \downarrow f \\
& & Y \end{array}
\]
where $\beta$ is a product of blowups of nonsingular subvarieties, so similarly we
have $$ \beta^*:H^i(Y, \co_Y) \cong H^i(W, \co_W) $$ for all $i$. This implies that
$f^{*}$ is one-to-one, and the theorem is proved.
\end{proof}

Suppose that $f: Y \rightarrow X$ is a morphism of schemes, and $\cf$ is a sheaf of
Abelian groups on $Y$. From the Leray spectral sequence
$H^i(X,R^jf_*\cf)\Rightarrow H^{i+j}(Y,\cf)$, we deduce the following exact sequence
\begin{equation}\label{eq9}
0\rightarrow H^1(X,f_*\cf)\rightarrow H^1(Y,\cf)\rightarrow
H^0(X,R^1f_*\cf)\rightarrow H^2(X,f_*\cf)\rightarrow H^2(Y,\cf).
\end{equation}

In the case of nonsingular varieties over a field of characteristic 0, we have a
good necessary and sufficient condition for the existence of an arithmetic
Macaulayfication.

\begin{thm} \label{nonsingular}
Suppose $X$ is a nonsingular projective scheme over a field $k$ of characteristic 0.
Then, $X$ has an arithmetic Macaulayfication if and only if
$H^0(X,\co_X)=k$ and $H^i(X, \co_X) = 0$ for
all $i =1, \ldots, \dim X-1$.
\end{thm}

\begin{proof} Suppose that $X$ is nonsingular and there exists an arithmetic
Macaulayfication $f:Y=\pr S\rightarrow X$. By Lemma \ref{cond}, we have
$H^0(Y,\co_Y)=k$ and $H^i(Y,\co_Y)=0$ for $0 < i < \dim Y$.

Let $g:Z\rightarrow Y$ be a resolution of singularites. Set $h=f\circ g$. Then,
$H^i(Z,\co_Z)\cong H^i(X,\co_X)$ for all $i$ by Theorem \ref{hiro}.

We have sequences
\[ H^i(X,\co_X)\stackrel{f^*}{\rightarrow}H^i(Y,\co_Y)\stackrel{g^*}{\rightarrow}
H^i(Z,\co_Z). \] $h^*=g^*\circ f^*$ is an isomorphism, so we have $H^0(X,\co_X)=k$
and $H^i(X,\co_X)=0$ for $0 < i < \dim X = \dim Y$. The necessary condition is
proved.

Now suppose that $H^0(X,\co_X)=k$ and $H^i(X, \co_X) = 0$ for  $i = 1, \ldots, \dim
X -1$. $X$ is a Cohen-Macaulay scheme since $X$ is nonsingular. Now using Lemma
\ref{cond}, we can embed $X$ as an arithmetically Cohen-Macaulay scheme $Y$. The
sufficient condition is proved.
\end{proof}

\begin{remark} The same proof shows that the conclusions of Theorem
\ref{nonsingular} hold if $X$ has rational singularities, over a field
of characteristic zero.
\end{remark}

From Kawasaki's work we easily deduce a very strong criterion for the existence of
an arithmetic Macaulayfication over a field of arbitrary characteristic.

\begin{thm}\label{ACM} Suppose that $X$ is an unmixed projective scheme of
dimension $\ge 1$ over a field $k$, $H^i(X,\co_X)=0$ for $1\le i\le\text{ dim }X-1$
and $H^0(X,\co_X)=k$. Then there exists an arithmetic Macaulayfication of $X$.
\end{thm}

\begin{proof}
$X=\pr R$ where $R=\oplus_{i\ge0}R_i$ is an unmixed, standard graded $k$-algebra.
Let $V$ be the (reduced) closed subscheme of $X$ of non Cohen-Macaulay points,
$s=\text{ dim }V$, $d=\text{ dim }X$, $z_1,\ldots,z_d$ be homogeneous elements of
$R$ satisfying the conclusions of Lemma 5.3  \cite{ka1}. Since $R$ is unmixed
$s<d-1$ (as follows from Corollary 2.4 \cite{ka}). Let $Q_i=(z_i,\ldots,z_d)\subset
R$ for $1\le i\le s+1$, and $I=Q_1\cdots Q_sQ_{s+1}^{d-s-1}\subset R$.

Suppose that $\alpha\in X$ is a closed point, $y\in R_1-\alpha$,
$x_i=\sfrac{z_i}{y^{\mbox{\small deg }z_i}}$ for $1\le i\le d$. Let
$q_i=(x_i,\ldots,x_d)$, $\beta=q_1\cdots q_sq_{s+1}^{d-s-1}\subset R_{(\alpha)}$. We
have $q_i=(Q_i)_{(\alpha)}$ and $\beta=I_{(\alpha)}$. If $\beta=R_{(\alpha)}$, then
$R_{(\alpha)}$ is Cohen-Macaulay and $R_{(\alpha)}[I_{(\alpha)}t]$ is
Cohen-Macaulay.  If $\beta\ne R_{(\alpha)}$, then there exists $l$ such that
$x_l,\ldots,x_d\in \alpha_{(\alpha)}$ and $x_{l-1}\not\in \alpha_{(\alpha)}$. As in
the proof of Theorem 5.1 \cite{ka1}, $x_l,\ldots, x_d$ is a subsystem of a
p-standard system of parameters for $R_{(\alpha)}$ and
$R_{(\alpha)}/(x_l,\ldots,x_d)R_{(\alpha)}$ is a Cohen-Macaulay ring if $l>1$.
$R_{(\alpha)}[\beta t]=R_{(\alpha)}[q_t\cdots q_sq_{s+1}^{d-s-1}t]$ is
Cohen-Macaulay by Corollary 4.5 \cite{ka}, since  $s<d-1$ and $((0):x_d)=(0)$ as
$R_{(\alpha)}$ is unmixed.

Let $\sh{I}$ be the sheafication of $I$, $Y=\pr (\oplus_{n\ge 0} \sh{I}^n)$, with
projection $\pi:Y\rightarrow \pr R=X$. For $\alpha\in X$ a closed point,
$R^i\pi_*\co_{Y,\alpha}=H^i(Y_{\alpha},\co_{Y_\alpha})$ where $Y_{\alpha}=\pr
R_{(\alpha)}[I_{(\alpha)}t] = Y \times_X \spec \co_{X, \alpha}$. Since
$R_{(\alpha)}[I_{(\alpha)}t]$ is Cohen-Macaulay, we have
$H^i(Y_{\alpha},\co_{Y_\alpha})=0$ for $i>0$ and
$(\pi_*\co_Y)_\alpha=R_{(\alpha)}=\co_{X,\alpha}$ by Theorem 4.1 \cite{l}. Thus
$R^i\pi_*\co _Y=0$ for $i>0$ and $\pi_*\co _Y=\co_X$.  From the Leray spectral
sequence we deduce that $$ H^i(Y,\co_Y)=H^i(X,\pi_*\co_Y)=H^i(X,\co_X)=0 $$ for
$1\le i\le d-1$ and $$ H^0(Y,\co_Y)=H^0(X,\pi_*\co_Y)=H^0(X,\co_X)=k. $$ It also
follows from what was shown above that if $\gamma \in Y$ is a closed point, $\co_{Y,
\gamma}$ is Cohen-Macaulay, so that $Y$ is a Cohen-Macaulay scheme. Lemma
\ref{cond} now implies that $Y=\pr S$ for some Cohen-Macaulay ring $S$.
\end{proof}

\begin{example}\label{Example1}
The converse of Theorem \ref{ACM} is not true, as can be seen from the following
simple example. Suppose that $k$ is an algebraically closed field, $X=\pr S$ is the
cuspidal plane curve with coordinate ring $S=k[y_0,y_1,y_2]/(y_0y_2^2-y_1^3)$.
$H^1(X,\co_X)\cong k$.  Let $Y=X\times\pp^1_k$. Note that $Y$ is a Cohen-Macaulay
scheme. $H^1(Y,\co_Y)\cong k\ne0$, by the K\"unneth formula. There is a natural
resolution of singularites $\pp^1_k\times \pp^1_k\rightarrow Y$, which is an
arithmetic Macaulayfication, as $\pp^1_k\times \pp^1_k\cong \pr R$, with
$R=k[x_0,x_1,x_2,x_3]/(x_0x_2-x_1x_3)$.
\end{example}

We observe that the converse of Theorem \ref{ACM} is true for normal projective
surfaces. For if $X$ is a projective normal surface and $f:Y\rightarrow X$ is an
arithmetic Macaulayfication, then $f_*\co_Y=\co_X$, so that
$k=H^0(Y,\co_Y)=H^0(X,\co_X)$, and $H^1(X,\co_X)=H^1(Y,\co_Y)=0$ by (\ref{eq9}) and
Lemma \ref{cond}.

The following example is of a normal 3-fold $X$ such that the converse of Theorem
\ref{ACM} is false.

\begin{example}\label{ExampleB}
There exists a normal projective 3-fold $B$ such that $H^2(B,\co_B)\ne 0$
and $B$ has an arithmetic Macaulayfication.
\end{example}

\begin{proof} In section III of \cite{C} an example is given of an $m$-primary ideal $I$
in the power series ring $\cc[[x,y,z]]$ such that $\oplus_{n\ge 0}I^n$ is normal
but not Cohen-Macaulay. The construction there yields an example of the desired
type.

Let $\beta :A\rightarrow \pp_\cc^3$ be the morphism obtained by blowing up a point
in $\pp_\cc^3$, and then blowing up the 12 points which are the intersection points
of a general hypersurface on the exceptional $\pp^2$ with a general cubic curve
$C''$ on the exceptional $\pp^2$. Let $C'$ be the strict transform of $C''$ on $A$.
In section III of \cite{C}, it is shown that there exists a projective morphism
$\alpha:A\rightarrow B$ such that $B$ is normal, $\alpha(C')$ is a point $Q$,
$A-C'\rightarrow B-Q$ is an isomorphism, $\alpha_*\co_A\cong\co_B$ and
$R^1\alpha_*\co_A\ne0$. Since $R^1\alpha_*\co_A$ is supported at the single point
$Q$, we have $H^0(B,R^1\alpha_*\co_A)\ne0$. Since $H^i(A,\co_A)\cong
H^i(\pp^3,\co_{\pp^3})$ for all $i$ (by Theorem \ref{hiro}), we have
$H^i(A,\co_A)=0$ for $i=1,2$ and $H^0(A,\co_A)=\cc$. Hence (by Lemma \ref{cond})
$A\cong \pr S$ where $S$ is a Cohen-Macaulay standard graded $\cc$-algebra. By
(\ref{eq9}), we have an isomorphism $H^2(B,\co_B)\cong H^0(B,R^1\alpha_*\co_A)\ne
0$.
\end{proof}

\section{Truncated Rees algebras}
Let $R$ be a standard graded $k$-algebra, $I \subseteq R$ a homogeneous ideal. The truncated
Rees algebras associated to $I$ are defined as follows.

\begin{mydef} Suppose that $I = \oplus_{t \ge \alpha}I_t$ is the homogeneous
decomposition of $I$, where $\alpha = \alpha(I)$ is the minimum degree in $I$. For
each $\lambda \ge \alpha$, we define the {\it truncated Rees algebra of $I$ at
degree $\lambda$} to be the Rees algebra
\[ \R_\lambda(I) = R[(I_\lambda) t] \subseteq R[t] \]
of the ideal generated by $I_\lambda$.
\end{mydef}

Define $\delta = \delta(I)$, the maximum degree of a minimal system of homogeneous
generators of $I$.

We will assume that $\lambda\ge\delta$. The truncated Rees algebra $\R_\lambda(I)$
has a  bi-gradation determined by $\deg F = (d,0)$ if $F\in R$ is homogeneous of
degree $d$, and $\deg t = (-\lambda,1)$, i.e.
\[ \R_\lambda(I)_{(p,q)} = I^q_{p+q\lambda}t^q. \]
It can be seen that
\[ R=\oplus_{n\ge 0}R_{\lambda}(I)_{(n,0)} \] as a graded subring of
$\R_\lambda(I)$, and
\[ S_{\lambda}=\oplus_{n\ge0}R_{\lambda}(I)_{(0,n)} \] is another
subring of $R_{\lambda}(I)$ which we will consider. There is a natural isomorphism
$S_{\lambda}\cong k[I_{\lambda}]$.

Set $X=\pr R$, $V_{\lambda}= \pr R_{\lambda}(I)$ (with respect to the above
bi-grading), $\overline V_{\lambda} =\pr S_{\lambda}$. We have canonical projections
$\pi_1:V_{\lambda}\rightarrow X$ and $\pi_2:V_{\lambda}
\rightarrow \overline V_{\lambda}$.

$V_{\lambda}$ can be  identified with the graph of the rational map $X \m \overline
V_{\lambda}$ induced by the natural inclusion $k[I_{\lambda}]\rightarrow R$, and we
have an isomorphism $V_{\lambda}\cong\pr (\oplus_{n\ge 0}\sh{I}^n)$, the blowup of
the sheafification $\sh{I}$ of $I$ (c.f. \cite{eh}). From now on we will assume that
$\lambda\ge\delta+1$. We then also have that $\overline V_{\lambda}\cong \pr
(\oplus_{n\ge 0}\sh{I}^n)$ (\cite[Lemma 1.1]{ch}) so that $V_{\lambda}\rightarrow
\overline V_{\lambda}$ is an isomorphism, and we have a natural diagram of
morphisms (where $\pi_2$ is an isomorphism): $$
\begin{array}{lllll}
&&V_{\lambda}&&\\
\pi_1&\swarrow&&\searrow&\pi_2\\
X&&\stackrel{\pi}{\leftarrow}&&\overline V_{\lambda}
\end{array}
$$

Let $\cl = \sh{I}\co_{\overline V_{\lambda}}$. The respective gradings on $R$, $S_{\lambda}$
and $R_{\lambda}(I)$ are related by  isomorphisms
$$
\cl^q\otimes\pi^*\co_X(q\lambda)\cong \co_{\overline V_{\lambda}}(q),
$$
\begin{equation}\label{eq*2}
\co_{V_{\lambda}}(p,q)\cong \pi_1^*\co_X(p)\otimes\pi_2^*\co_{\overline V_{\lambda}}(q).
\end{equation}

For $q\in\zz$, let $\cm_q$ be the sheafification on $X$ of the graded $R$-module $$
M_q=\oplus_{i\ge 0}R_{\lambda}(I)_{(i,q)} $$ so that (since $\lambda\ge\delta$) $$
\cm_q=\left\{\begin{array}{ll}
\co_X&\text{if }q=0\\
\sh{I}^q(\lambda q)&q>0\\
0&q<0
\end{array}\right.
$$
Thus for $p\in \zz$,
$$
\cm_q(p)=\left\{\begin{array}{ll}
\co_X(p)&\text{if }q=0\\
\sh{I}^q(\lambda q+p)&q>0\\
0&q<0
\end{array}\right.
$$

For $p\in\zz$, let $\cn_p$ be the sheafification on $\overline V_\lambda$ of the
graded $S_\lambda$ module $$ N_p=\oplus_{i\ge 0}R_{\lambda}(I)_{(p,i)}. $$ Observe
that $N_0=k[I_\lambda]=S_\lambda$, $$
\cn_p=\left\{\begin{array}{ll}
\co_{\overline V_\lambda}&\text{if }p=0\\
\pi^*\co_X(p)&p>0\\
0&p<0
\end{array}\right.
$$
Thus for $q\in \zz$,
$$
\cn_p(q)=\left\{\begin{array}{ll}
\pi^*\co_X(p)\otimes \co_{\overline V_\lambda}(q)&\text{if }p\ge0\\
0&p<0\\
\end{array}\right.
$$

Our main result in this section is to show that for a certain class of standard
graded $k$-algebras $R$ and homogeneous ideals $I$, the truncated Rees algebras
$\R_\lambda(I)$ of $I$ are  Cohen-Macaulay for large $\lambda$.

\begin{thm} \label{trunc-rees}
Suppose that $R$ is a Cohen-Macaulay standard graded $k$-algebra of positive
dimension $d+1$ that has negative $a$-invariance, $a(R) < 0$ (Since $R$ is
Cohen-Macaulay this is equivalent to $H^{d+1}_{{\goth{m}_1}}(R)_p = 0 \ \forall p
\ge 0$, where $\goth{m}_1$ is the maximal ideal of $R$). Let $I \subseteq R$ be a
homogeneous ideal, and suppose that $\lambda\ge\delta(I)+1$.

Let $\sh{I}$ be the ideal sheaf associated to $I$ on $X=\pr R$,
$$
E\cong \pr \oplus_{n\ge 0}
\sh{I}^n/\sh{I}^{n+1}
$$
 be the exceptional divisor of $\pi_1:V_{\lambda}\rightarrow X$, with dualizing sheaf
$\omega_E$ on $E$. Suppose that
\begin{equation}\label{eq6}
\left\{ \begin{array}{rcl}
\pi_{1*} \co_E(-\lambda m,m) & = & \sh{I}^m/\sh{I}^{m+1}, \forall \ m \ge 0, \\
R^i \pi_{1*} \co_E(-\lambda m,m) & = & 0, \forall \ i > 0, m \ge 0, \\
R^i \pi_{1*} \omega_E(-\lambda m,m) & = & 0, \forall \ i > 0, m \ge 2
\end{array} \right.
\end{equation}
Then, there exists an integer $\lambda_0$ such that for all
$$
\lambda \ge \lambda_0\ge\delta+1
$$
the truncated Rees algebra $\R_\lambda(I)$ is Cohen-Macaulay.
\end{thm}

To prove Theorem \ref{trunc-rees}, we shall combine the method of \cite{hyry} for
studying the local cohomology of multi-graded algebras with the results of
\cite{ch}. Suppose that $\lambda\ge\delta+1$. For convenience, denote
$S_{V_\lambda} = \R_\lambda(I)$. We need to show that $S_{V_\lambda}$ is a
Cohen-Macaulay ring for $\lambda \gg 0$.

Let
$$
\goth{m}_1=\oplus_{i>0}R_{\lambda}(I)_{(i,0)}
$$
 be the irrelevant ideal of $R$,
$$
\goth{n}_1=\goth{m}_1R_{\lambda}(I)=\oplus_{i>0,j\ge 0}R_{\lambda}(I)_{(i,j)}.
$$ Let $$
\goth{m}_2=\oplus_{j>0}R_{\lambda}(I)_{(0,j)}
$$
 be the irrelevant ideal of $S_{\lambda}$,
$$
\goth{n}_2=\goth{m}_2R_{\lambda}(I)=\oplus_{i\ge0,j>0}R_{\lambda}(I)_{(i,j)}.
$$ Let $$
\goth{m}=\oplus_{i+j>0}R_{\lambda}(I)_{(i,j)}
$$
 and
$$
\goth{n}=\oplus_{i,j>0}R_{\lambda}(I)_{(i,j)}.
$$ Then $\goth{n}_1+\goth{n}_2=\goth{m}$ and $\goth{n}_1\cap\goth{n}_2=\goth{n}$.

$S_{V_\lambda}$ is Cohen-Macaulay if and only if
\[ H^i_{\goth{m}}(S_{V_\lambda}) = 0, \ \forall \ i = 0, \ldots, d+1, \]
where $d = \dim X = \dim S_{V_\lambda} - 2$,
so that from the Mayer-Vietoris sequence of cohomologies,
\[ \ldots \rightarrow H^i_\goth{m}(S_{V_\lambda}) \rightarrow H^i_{\goth{n}_1}(S_{V_\lambda})
\oplus H^i_{\goth{n}_2}(S_{V_\lambda}) \rightarrow H^i_\goth{n}(S_{V_\lambda}) \rightarrow
H^{i+1}_\goth{m}(S_{V_\lambda}) \rightarrow \ldots. \]
we see that
$S_{V_{\lambda}}$ is Cohen-Macaulay if and only if
\begin{equation}\label{eq0}
\left\{ \begin{array}{rclcl}
H^i_{\goth{n}_1}(S_{V_\lambda}) & \oplus & H^i_{\goth{n}_2}(S_{V_\lambda}) &
\stackrel{\sim}{\rightarrow} &
H^i_\goth{n}(S_{V_\lambda}), \ \forall \ i = 0, \ldots, d \\
H^{d+1}_{\goth{n}_1}(S_{V_\lambda}) & \oplus & H^{d+1}_{\goth{n}_2}(S_{V_\lambda}) &
\hookrightarrow & H^{d+1}_\goth{n}(S_{V_\lambda})
\end{array} \right.
\end{equation}

We have isomorphisms
$$
H^i_{\goth{m}_1}(M_q)_p\cong H^i_{\goth{n}_1}(S_{V_{\lambda}})_{(p,q)}
$$
and
$$
H^i_{\goth{m}_2}(N_p)_q\cong H^i_{\goth{n}_2}(S_{V_{\lambda}})_{(p,q)}
$$
for all $p,q\in\zz$ (c.f. Lemma 2.1 \cite{cht}).

For $p,q\in\zz$, we have commutative diagrams with exact rows \cite[Theorem
1.4]{hyry}.
\begin{equation}\label{eq1}
\begin{array}{ccccccccc}
0 \rightarrow & H^0_\goth{n}(S_{V_\lambda})_{(p,q)} & \rightarrow &
(S_{V_\lambda})_{(p,q)} & \rightarrow & H^0({V_\lambda}, \co_{V_\lambda}(p,q)) &
\rightarrow & H^1_\goth{n}(S_{V_\lambda})_{(p,q)}
& \rightarrow 0 \\
& \uparrow & & \uparrow\wr & & \uparrow & & \uparrow  & \\
0 \rightarrow & H^0_{\goth{n}_1}(S_{V_\lambda})_{(p,q)} & \rightarrow &
(S_{V_\lambda})_{(p,q)} & \rightarrow & H^0(X, \cm_q(p)) &
\rightarrow & H^1_{\goth{n}_1}(S_{V_\lambda})_{(p,q)} & \rightarrow 0 \end{array}
\end{equation}
and isomorphisms
\begin{equation}\label{eq3}
H^i_{\goth{n}}(S_{V_{\lambda}})_{(p,q)}\cong H^{i-1}(V_{\lambda},\co_{V_{\lambda}}(p,q))
\end{equation}
and
\begin{equation}\label{eq4}
H^i_{\goth{n}_1}(S_{V_{\lambda}})_{(p,q)}\cong H^{i-1}(X,\cm_q(p))
\end{equation}
for all $i\ge 2$.

For $p,q\in\zz$, we have commutative diagrams with exact rows (\cite[Theorem
1.4]{hyry}).
\begin{equation}\label{eq2}
\begin{array}{ccccccccc}
0 \rightarrow & H^0_\goth{n}(S_{V_\lambda})_{(p,q)} & \rightarrow &
(S_{V_\lambda})_{(p,q)} & \rightarrow & H^0({V_\lambda}, \co_{V_\lambda}(p,q)) &
\rightarrow & H^1_\goth{n}(S_{V_\lambda})_{(p,q)}
& \rightarrow 0 \\
&  \uparrow & & \uparrow\wr & & \uparrow & & \uparrow  & \\
0 \rightarrow & H^0_{\goth{n}_2}(S_{V_\lambda})_{(p,q)} & \rightarrow &
(S_{V_\lambda})_{(p,q)} & \rightarrow & H^0(\overline V_{\lambda}, \cn_p(q)) &
\rightarrow & H^1_{\goth{n}_1}(S_{V_\lambda})_{(p,q)} & \rightarrow 0 \end{array}
\end{equation}
and isomorphisms
\begin{equation}\label{eq5}
H^i_{\goth{n}_2}(S_{V_{\lambda}})_{(p,q)}\cong H^{i-1}(\overline V_{\lambda},\cn_p(q))
\end{equation}
for all $i\ge 2$.

By Lemma 2.1 \cite{ch}
$$
R^i\pi_{1*}\co_{V_\lambda}(p,q)=0\text{ for }i>0,q\ge 0,p\in\zz
$$
and
$$
\pi_{1*}\co_{V_\lambda}(p,q)\cong\sh{I}^q(p+q\lambda)\text{ for }q\ge 0, p\in\zz.
$$
Thus by the Leray spectral sequence,
\begin{equation}\label{eq*1}
H^i(\overline V_\lambda,\co_{\overline V_\lambda}(q)\otimes\pi^*\co_X(p))\cong
H^i(V_{\lambda},\co_{V_{\lambda}}(p,q))\cong H^i(X,\sh{I}^q(p+q\lambda))
\end{equation}
for $i\ge 0, q\ge 0,p\in\zz$.

\begin{pro}\label{TheoremBG}There exists
$\lambda_0\ge\delta+1$ such that for $\lambda\ge \lambda_0$,
\begin{enumerate}
\item $H^0(V_{\lambda},\co_{V_{\lambda}}(p,q))=(S_{V_{\lambda}})_{(p,q)}$
for $p,q\ge 0$.
\item $H^0(V_{\lambda},\co_{V_{\lambda}}(p,q))=0$ for $p,q<0$
\item $H^i(V_{\lambda},\co_{V_{\lambda}}(p,q))=0$ for $i>0,p\ge0,q>0$.
\item $H^i(V_{\lambda},\co_{V_{\lambda}}(p,q))=0$ for $i<d,p,q<0$
\end{enumerate}
\end{pro}

\begin{proof} (1) follows from (\ref{eq*2}) and Lemma 1.3 \cite{ch}.

We now prove (2). After possibly tensoring with an extension field of $k$, we may
suppose that $k$ is an infinite field. Suppose that $\lambda>\delta+1$. Then
$\co_{V_\lambda}(-r(\lambda-(\delta+1)),r)$ is very ample for $r\ge 1$, by
(\ref{eq*2}) and Lemma 1.1 \cite{ch}. For $r\ge 1$, there exists $\sigma_r\in
H^0(V_\lambda,\co_{V_\lambda}(-r(\lambda-(\delta+1)),r))$ such that $(\sigma_r)$
contains no associated  primes of $V_\lambda$. Thus $$
\sigma_r:\co_{V_\lambda}(r(\lambda-(\delta+1)),-r)\rightarrow \co_{V_\lambda}
$$
are inclusions. For $s>0$, we have inclusions
$$
\co_{V_\lambda}(r(\lambda-(\delta+1))-s,-r)\rightarrow \co_{V_\lambda}(-s,0)
$$
which induce inclusions
$$
\pi_{1*}\co_{V_\lambda}(r(\lambda-(d+1))-s,-r)\rightarrow \co_X(-s)
$$ and thus  inclusions $$
H^0(V_{\lambda},\co_{V_\lambda}(r(\lambda-(\delta+1))-s,-r))\rightarrow
H^0(X,\co_X(-s)) $$ $H^0(X,\co_X(-s))=0$ for $s<0$ since $R$ is Cohen-Macaulay.
Since $\lambda>\delta+1$, we have $$ H^0(V_\lambda,\co_{V_\lambda}(p,q))=0 $$ if
$p,q<0$.

(3) follows from (\ref{eq*2}) and Proposition 3.1 \cite{ch}.

(4) follows from (\ref{eq*2}), Lemma 2.2 \cite{ch} and Proposition 3.2 \cite{ch}.
\end{proof}

We also have
\begin{equation}\label{eq7}
H^i(X,\co_X(p))=0 \text{ for }i>0,p\ge 0
\end{equation}
since $a(R)<0$, and
\begin{equation}\label{eq8}
H^i_{\goth{m}_1}(R)=0 \text{ for }i<d+1,
\end{equation}
since $R$ is Cohen-Macaulay, so that $$ H^0(X,\co_X(p))=\left\{\begin{array}{ll}
0&p<0\\
(S_{V_\lambda})_{(p,0)}&p\ge 0.
\end{array}\right.
$$

The conditions of (\ref{eq0}) now follow for $\lambda\ge\lambda_0$ by (\ref{eq*1}),
(\ref{eq*2}), Proposition \ref{TheoremBG}, (\ref{eq7}), (\ref{eq8}) and (\ref{eq1})
- (\ref{eq5}).  We have thus finished the proof Theorem \ref{trunc-rees}.

\begin{cor}\label{CorollaryLCI} Suppose that $X$ is a projective
Cohen-Macaulay scheme over a field $k$ such that $H^0(X,\co_X)=k$,
$H^i(X,\co_X)=0$ for $i>0$ and $\sh{I}\subset \co_X$ is an ideal sheaf which is locally
a complete intersection. Then
\begin{enumerate}
\item There exists a Cohen-Macaulay standard graded $k$-algebra $R$ with $a(R)<0$
such that $X\cong \pr R$.
\item If $I\subset R$ is a homogeneous ideal such that $\tilde I\cong \sh{I}$,
then there exists $\lambda_0\ge\delta(I)+1$ (where $\delta(I)$ is the maximum degree
of a minimal set of homogeneous generators of $I$) such that the truncated Rees
algebra
$$
\R_\lambda(I) = R[(I_\lambda) t]
$$
is Cohen-Macaulay for $\lambda\ge\lambda_0$.
\end{enumerate}
\end{cor}

\begin{proof} (1) follows from Lemma \ref{cond}. $I$ satisfies the condition
(\ref{eq6}) by (\ref{eq*2}) and Example 2.3 of \cite{ch}, so that (2) follows from
Theorem \ref{trunc-rees}.
\end{proof}

$X=\pp_k^n$ is an especially important example satisfying the conditions
on $X$ of the corollary.
Other classes of situations where the conditions in (\ref{eq6}) are satisfied
can be found from \cite{ch}.

\begin{remark} Since $R$ is Cohen-Macaulay, the hypothesis $a(R) < 0$ in Theorem
\ref{trunc-rees} is equivalent to $H^{d+1}_{{\goth{m}_1}}(R)_p = 0 \ \forall p
\ge 0$. This is an extra condition compared to \cite[Theorem 4.1]{ch}.
\end{remark}

\begin{remark} Suppose that $R$ is Cohen-Macaulay, $I\subset R$ is homogeneous, and
(\ref{eq6}) holds. If there exists $\lambda\ge\delta(I)+1$ such that
$\R_\lambda(I)$ is Cohen-Macaulay, then $a(R)<0$.
\end{remark}

\begin{proof} Let notation be as in the proof of Theorem \ref{trunc-rees}. The
Remark follows from (\ref{eq0}), and the fact that for $p\ge 0$,
$$
H_{\goth{n}_1}^{d+1}(S_{V_{\lambda}})_{(p,0)}\cong H^d(X,\co_X(p)),
$$
$$
\begin{array}{ll}
H_{\goth{n}_2}^{d+1}(S_{V_\lambda})_{(p,0)}&\cong H^d(\overline V_\lambda,\pi^*\co_X(p))\\
&\cong H^d(X,\co_X(p))
\end{array}
$$
and
$$
H^{d+1}_{\goth{n}}(S_{V_\lambda})_{(p,0)}\cong H^d(V_\lambda,\co_{V_\lambda}(p,0))\cong
H^d(X,\co_X(p)).
$$
\end{proof}

\end{document}